\documentclass[12pt]{article}
\usepackage{amsmath,amssymb,theorem}
\usepackage{graphicx}
\usepackage{caption}

\usepackage{subfigure}
\usepackage{enumerate}

\usepackage{psfrag}
\usepackage{mathrsfs}
\usepackage{tikz-cd}
\usepackage{bbm}

\usepackage{epstopdf}
\usepackage{float}

\theorembodyfont{\rmfamily}
\newtheorem{thm}{Theorem}[section]
\newtheorem{conj}[thm]{Conjecture}
\newtheorem{cor}[thm]{Corollary}
\newtheorem{lem}[thm]{Lemma}
\newtheorem{prop}[thm]{Proposition}
\theorembodyfont{\rmfamily}

\def\pf{\bigskip\noindent {\bf Proof.}~~}

\def\mytextindent#1{\indent\llap{#1\enspace}\ignorespaces}
\def\fattextindent#1{\indent\indent\llap{#1\enspace}\ignorespaces}
\def\myitem{\par\hangindent\parindent\mytextindent}
\def\myitemitem{\par\hangindent\parindent\fattextindent}
\def\proofsquare{ \hfill\vrule height3pt width6pt depth2pt}

\def\dfn#1{{\sl #1}}

\def\less{\backslash}

\newcounter{counter}

\begin{document}
\title{Double-critical   graph conjecture for claw-free graphs}
\author{Martin Rolek\thanks{E-mail address: mrolek@knights.ucf.edu.}
~and~ Zi-Xia Song\thanks{Corresponding author. E-mail address: Zixia.Song@ucf.edu.}\\
Department of Mathematics\\
University of Central Florida\\
 Orlando, FL 32816
}

\maketitle
\begin{abstract}
A connected graph $G$ with chromatic number $t$ is {\sl double-critical} if $G \backslash \{x, y\}$ is $(t - 2)$-colorable for each edge $xy \in E(G)$.
The complete graphs are the only known examples of double-critical graphs.
A long-standing conjecture of Erd\H os and Lov\'asz from 1966, which is referred to as the {\sl Double-Critical Graph Conjecture}, states that there are no other double-critical graphs.
That is, if a graph $G$ with chromatic number $t$ is double-critical, then $G$ is the complete graph on $t$ vertices.
This has been verified for $t \le 5$, but remains open for $t \ge 6$.
In this paper, we first prove that if $G$ is a non-complete, double-critical graph with chromatic number $t \ge 6$, then no vertex of degree $t + 1$ is adjacent to a vertex of degree $t+1$, $t + 2$, or $t + 3$ in $G$.
We then use this result to show that the Double-Critical Graph Conjecture is true for double-critical graphs $G$ with chromatic number $t \le 8$ if $G$ is claw-free.
\end{abstract}


\baselineskip 18pt

\section{Introduction}
All graphs considered in this paper are finite and without loops or multiple edges.
For a graph $G$, we will use $V(G)$ to denote the vertex set, $E(G)$ the edge set, $e(G)$ the number of edges, $\alpha(G)$ the independence number, $\omega(G)$ the clique number, $\chi(G)$ the chromatic number, and $\overline{G}$ the complement of $G$, respectively.
For a vertex $x \in V(G)$, we will use $N_G(x)$ to denote the set of vertices in $G$ which are adjacent to $x$.
We define $N_G[x] = N_G(x) \cup \{x\}$ and $d_G(x) = |N_G(x)|$.
Given vertex sets $A, B \subseteq V(G)$, we say that $A$ is \dfn{complete to} (resp. \dfn{anti-complete to}) $B$ if for every $a \in A$ and every $b \in B$, $ab \in E(G)$ (resp. $ab \notin E(G)$).
The subgraph of $G$ induced by $A$, denoted $G[A]$, is the graph with vertex set $A$ and edge set $\{xy \in E(G) : x, y \in A\}$. We denote by $B \less A$ the set $B - A$, $e_G(A, B)$ the number of edges between $A$ and $B$ in $G$, and $G \less A$ the subgraph of $G$ induced on $V(G) \less A$, respectively.
If $A = \{a\}$, we simply write $B \less a$, $e_G(a, B)$, and $G \less a$, respectively.  A graph $H$ is an \dfn{induced subgraph} of a graph $G$ if $V(H)\subseteq V(G)$   and $H=G[V(H)]$.  A graph $G$ is \dfn{claw-free} if $G$ does not contain $K_{1,3}$ as an induced subgraph. 
Given two graphs $G$ and $H$, the \dfn{union} of $G$ and $H$, denoted $G \cup H$, is the graph with vertex set $V(G) \cup V(H)$ and edge set $E(G) \cup E(H)$.
Given two isomorphic graphs $G$ and $H$, we may (with a slight but common abuse of notation) write $G = H$.    A cycle with $t \ge 3$ vertices is denoted by $C_t$. Throughout this paper, a proper vertex coloring of a graph $G$ with $k$ colors is called a \dfn{$k$-coloring of $G$}.
\medskip

In 1966, the following conjecture of Lov\'asz was published by Erd\H os~\cite{Erdos1968} and is known as the Erd\H os-Lov\'asz Tihany Conjecture.

\begin{conj}\label{conj:ELTC}
For any integers $s, t \ge 2$ and any graph $G$ with $\omega(G) < \chi(G) = s + t - 1$, there exist disjoint subgraphs $G_1$ and $G_2$ of $G$ such that $\chi(G_1) \ge s$ and $\chi(G_2) \ge t$.
\end{conj}

To date, Conjecture~\ref{conj:ELTC} has been shown to be true only for values of $(s, t) \in \{ (2, 2), \allowbreak (2, 3), \allowbreak (2, 4), \allowbreak (3, 3), \allowbreak (3, 4), \allowbreak (3, 5) \}$.
The case $(2, 2)$ is trivial.
The case $(2, 3)$ was shown by Brown and Jung in 1969~\cite{Brown1969}.
Mozhan~\cite{Mozhan1987} and Stiebitz~\cite{Stiebitz1987a} each independently showed the case $(2, 4)$ in 1987.
The cases $(3, 3)$, $(3, 4)$, and $(3, 5)$ were also settled by Stiebitz in 1987 \cite{Stiebitz1987b}.
Recent work on the Erd\H os-Lov\'asz Tihany Conjecture has focused on proving the conjecture for certain classes of graphs.
Kostochka and Stiebitz~\cite{Kostochka2008} showed the conjecture holds for line graphs.
Balogh, Kostochka, Prince, and Stiebitz~\cite{Balogh2009} then showed that the conjecture holds for all quasi-line graphs and all graphs $G$ with $\alpha(G) = 2$.
More recently, Chudnovsky, Fradkin, and Plumettaz~\cite{Chudnovsky2013} proved the following slight weaking of Conjecture~\ref{conj:ELTC} for claw-free graphs, the proof of which is long and relies heavily on the structure theorem for claw-free graphs developed by Chudnovsky and Seymour~\cite{Chudnovsky2008}.

\begin{thm}\label{thm:ELTCClawFree}
Let $G$ be a claw-free graph with $\chi(G) > \omega(G)$.
Then there exists a clique $K$ with $|V(K)| \le 5$ such that $\chi(G \setminus V(K)) > \chi(G) - |V(K)|$.
\end{thm}

The most recent result related to the Erd\H os-Lov\'asz Tihany Conjecture is due to Stiebitz~\cite{Stiebitz2016}, who showed that for integers $s, t \ge 2$, any graph $G$ with $\omega(G) < \chi(G) = s + t - 1$ contains disjoint subgraphs $G_1$ and $G_2$ of $G$ with either $\chi(G_1) \ge s$ and $\text{col}(G_2) \ge t$, or $\text{col}(G_1) \ge s$ and $\chi(G_2) \ge t$, where $\text{col}(H)$ denotes the coloring number of a graph $H$.
\medskip

If we restrict $s = 2$ in Conjecture~\ref{conj:ELTC}, then the Erd\H os-Lov\'asz Tihany Conjecture states that for any graph $G$ with $\chi(G) > \omega(G) \ge 2$, there exists an edge $xy \in E(G)$ such that $\chi(G \less \{x, y\}) \ge \chi(G) - 1$.
To prove this special case of Conjecture~\ref{conj:ELTC}, suppose for a contradiction that no such edge exists.
Then $\chi(G \less \{x, y\}) = \chi(G) - 2$ for every edge $xy \in E(G)$.
This motivates the definition of double-critical graphs.
A connected graph $G$ is \dfn{double-critical} if for every edge $xy \in E(G)$, $\chi(G \less \{x, y\}) = \chi(G) - 2$.
A graph $G$ is \dfn{$t$-chromatic} if $\chi(G) = t$.
We are now ready to state the following conjecture, which is referred to as the \dfn{Double-Critical Graph Conjecture}, due to Erd\H os and Lov\'asz~\cite{Erdos1968}.

\begin{conj}\label{conj:DC}
Let $G$ be a double-critical, $t$-chromatic graph.
Then $G = K_t$.
\end{conj}

Since Conjecture~\ref{conj:DC} is a special case of Conjecture~\ref{conj:ELTC}, it has been settled in the affirmative for $t \le 5$~\cite{Mozhan1987, Stiebitz1987a}, for line graphs~\cite{Kostochka2008}, and for quasi-line graphs and graphs with independence number two~\cite{Balogh2009}.
Representing a weakening of Conjecture~\ref{conj:DC}, Kawarabayashi, Pedersen, and Toft~\cite{Kawarabayashi2010} have shown that any double-critical, $t$-chromatic graph contains $K_t$ as a minor for $t \in \{6, 7\}$.
As a further weakening, Pedersen~\cite{Pedersen2011} showed that any double-critical, $8$-chromatic graph contains $K_8^-$ as a minor.
Albar and Gon\c calves~\cite{Albar2015} later proved that any double-critical, 8-chromatic graph contains $K_8$ as a minor.
Their proof is computer-assisted.
The present authors~\cite{Rolek2016a} gave a computer-free proof of the same result and further showed that any double-critical, $t$-chromatic graph contains $K_9$ as a minor for all $t \ge 9$.
We note here that Theorem~\ref{thm:ELTCClawFree} does not completely settle Conjecture~\ref{conj:DC} for all claw-free graphs.
Recently, Huang and Yu~\cite{Huang2016} proved that the only double-critical, $6$-chromatic, claw-free graph is $K_6$.
We prove the following main results in this paper.
Theorem~\ref{thm:t+1nbrs} is a generalization of a result obtained in \cite{Kawarabayashi2010} that no two vertices of degree $t+1$ are adjacent in any non-complete, double-critical, $t$-chromatic graph.

\begin{thm}\label{thm:t+1nbrs}
If $G$ is a non-complete, double-critical, $t$-chromatic graph with $t \ge 6$, then for any vertex $x \in V(G)$ with $d_G(x) = t + 1$, the following hold:
\medskip

\myitem{(a)} $e(\overline{G[N_G(x)]}) \ge 8$; and
\smallskip

\myitem{(b)} for any vertex $y \in N_G(x)$, $d_G(y) \ge t + 4$.
Furthermore, if $d_G(y) = t + 4$,  then $|N_G(x) \cap N_G(y)| = t-2$ and   $\overline{G[N_G(x)]}$ contains   either only  one cycle, which is isomorphic to  $C_8$,  or  exactly  two cycles, each of which is isomorphic to   $C_{5}$. 
\end{thm}
\medskip

\noindent Corollary~\ref{cor:t+1nbrs} below follows immediately from Theorem~\ref{thm:t+1nbrs}.
 
\begin{cor}\label{cor:t+1nbrs}
If $G$ is a non-complete, double-critical, $t$-chromatic graph with $t \ge 6$, then no vertex of degree $t + 1$ is adjacent to a vertex of degree $t + 1$, $t + 2$, or $t + 3$ in $G$.
\end{cor}

We then use Corollary~\ref{cor:t+1nbrs} to prove the following main result.

\begin{thm}\label{thm:DCClawFree678}
Let $G$ be a double-critical, $t$-chromatic graph with $t \in \{6, 7, 8\}$.
If $G$ is claw-free, then $G = K_t$.
\end{thm}
\medskip

The rest of this paper is organized as follows.
In Section~\ref{sec:Prelim}, we first list some known properties of non-complete, double-critical graphs obtained in~\cite{Kawarabayashi2010} and then establish a few new ones.
In particular, Lemma~\ref{lem:KempeNbrs} turns out to be very useful.
Our new lemmas lead to a very short proof of Theorem~\ref{thm:DCClawFree678} for $t = 6, 7$, which we place at the end of Section~\ref{sec:Prelim}.
We prove the remainder of our main results in Section~\ref{sec:DCClawFree678}.

\section{Preliminaries}\label{sec:Prelim}

The following is a summary of the basic properties of non-complete, double-critical graphs shown by Kawarabayashi, Pedersen, and Toft in \cite{Kawarabayashi2010}.
\medskip

\begin{prop}\label{prop:DCBasics}
If  $G$ is a non-complete, double-critical, $t$-chromatic graph, then all of the following are true.

\myitem{(a)} $G$ does not contain $K_{t - 1}$ as a subgraph.
\smallskip

\myitem{(b)} For all edges $xy$, every $(t - 2)$-coloring $c: V(G) \less \{x, y\} \rightarrow \{1, 2, \dots, t - 2\}$ of $G \less \{x, y\}$, and any non-empty sequence $j_1, j_2, \dots, j_i$ of $i$ different colors from $\{1, 2, \dots, t - 2\}$, there is a path of order $i + 2$ with vertices $x, v_1, v_2, \dots, v_i, y$ in order such that $c(v_k) = j_k$ for all $k \in \{1, 2, \dots, i\}$.
\smallskip

\myitem{(c)} For any edge $xy \in E(G)$, $x$ and $y$ have at least one common neighbor in every color class of any $(t - 2)$-coloring of $G \less \{x, y\}$.
In particular, every edge $xy \in E(G)$ belongs to at least $t - 2$ triangles.
\smallskip

\myitem{(d)} There exists at least one edge $xy \in E(G)$ such that $x$ and $y$ share a common non-neighbor in $G$.
\smallskip

\myitem{(e)} For any edge $xy \in E(G)$,  the subgraph of $G$ induced by $N_G(x) \less N_G[y]$ contains no isolated vertices.
 In particular, no vertex  of $N_G(x)$ can have degree one in $\overline{G[N_G(x)]}$.
\smallskip

\myitem{(f)} $\delta(G) \ge t + 1$.
\smallskip

\myitem{(g)} For any vertex $x \in V(G)$, $\alpha(G[N_G(x)]) \le d_G(x) - t + 1$.
\smallskip

\myitem{(h)} For any vertex $x$ with at least one non-neighbor in $G$, $\chi(G[N_G(x)]) \le t - 3$.
\smallskip

\myitem{(i)} For any $x \in V(G)$ with $d_G(x) = t + 1$, $\overline{G[N_G(x)]}$ is the union of isolated vertices and cycles of length at least five.
Furthermore, there must be at least one such cycle in $\overline{G[N_G(x)]}$.
\smallskip

\myitem{(j)} No two vertices of degree $t + 1$ are adjacent in $G$.
\end{prop}

\medskip
\noindent
We next establish some new properties of non-complete, double-critical graphs. 
 
\begin{lem}\label{lem:DomVertex}
Let $G$ be a double-critical, $t$-chromatic graph and let $x \in V(G)$.
If $d_G(x) = |V(G)| - 1$, then $G \less x$ is a double-critical, $(t - 1)$-chromatic graph.
\end{lem}

\pf  Let $uv$ be any edge of $G \less x$.
Clearly, $\chi(G \less x) = t - 1$.
Since $G$ is double-critical, $\chi(G \less \{u, v\}) = t - 2$ and so $\chi(G \less \{u, v, x\}) = t - 3$ because $x$ is adjacent to all the other vertices in $G \less \{u, v\}$.
Hence $G \less x$ is double-critical and $(t - 1)$-chromatic.
\proofsquare\medskip

\begin{lem}\label{lem:ColorNonNbrs}
If $G$ is a non-complete, double-critical, $t$-chromatic graph, then for any $x \in V(G)$ with at least one non-neighbor in $G$, $\chi(G \less N_G[x]) \ge 3$.
In particular, $G \less N_G[x]$ must contain an odd cycle, and so $d_G(x) \le |V(G)| - 4$.
\end{lem}

\pf Let $x$ be any vertex in $G$ with $d_G(x) < |V(G)| - 1$ and let $H = G \less N_G[x]$.
Suppose that $\chi(H) \le 2$.
Since $d_G(x) < |V(G)| - 1$, $H$ contains at least one vertex.
Let $y \in V(H)$ be adjacent to a vertex $z \in N_G(x)$.
This is possible because $G$ is connected. 
If $H$ has no edge, then $G \less (V(H) \cup \{z\})$ has a $(t - 2)$-coloring $c$, which can be extended to a $(t - 1)$-coloring of $G$ by assigning all vertices in $V(H)$ the color $c(x)$ and assigning a new color to the vertex $z$, a contradiction.
Thus $H$ must contain at least one edge, and so $\chi(H) = 2$.
Let $(A, B)$ be a bipartition of $H$.
Now $G \less H$ has a $(t - 2)$-coloring $c^*$, which again can be extended to a $(t - 1)$-coloring of $G$ by assigning all vertices in $A$ the color $c^*(x)$ and all vertices in $B$ the same new color, a contradiction.
This proves that $\chi(H) \ge 3$, and so $H$ must contain an odd cycle.
Therefore $d_G(x) \le |V(G)| - 4$.
\proofsquare
\medskip

\begin{lem}\label{lem:KempeNbrs}
Let $G$ be a double-critical, $t$-chromatic graph.
For any edge $xy \in E(G)$, let $c$ be any $(t - 2)$-coloring of $G \less \{x, y\}$ with color classes $V_1, V_2, \dots, V_{t - 2}$.
Then the following two statements are true.
\medskip

\myitem{(a)} For any $i, j \in \{1, 2, \dots, t - 2\}$ with $i \ne j$, if $N_G(x) \cap N_G(y) \cap V_i$ is anti-complete to $N_G(x) \cap V_j$, then there exists at least one edge between $(N_G(y) \less N_G(x)) \cap V_i$ and $N_G(x) \cap V_j$ in $G$.
In particular, $(N_G(y) \less N_G(x)) \cap V_i \ne \emptyset$.
\medskip

\myitem{(b)} Assume that  $d_G(x) = t + 1$ and $y$ belongs to a cycle of length $k \ge 5$ in $\overline{G[N_G(x)]}$.

\myitemitem{($b_1$)} If $k \ge 7$, then $d_G(y) \ge t + e(\overline{G[N(x)]}) - 4$;

\myitemitem{($b_2$)}  If $k = 6$, then $d_G(y) \ge \max\{t + 2, t + e(\overline{G[N_G(x)]}) - 5\}$; and

\myitemitem{($b_3$)}  If $k = 5$, then $d_G(y) \ge \max\{t + 2, t + e(\overline{G[N_G(x)]}) - 6\}$.
\end{lem}

\pf Let $G, x, y, c$ be as given in the statement.  For any $i, j \in \{1, 2, \dots, t - 2\}$ with $i \ne j$, if $N_G(x) \cap N_G(y) \cap V_i$ is anti-complete to $N_G(x) \cap V_j$, 
then $G$ is non-complete. 
By Proposition~\ref{prop:DCBasics}(b), there must exist a path $x, u_j, u_i, y$ in $G$ such that $c(u_j) = j$ and $c(u_i) = i$.
Clearly, $u_j u_i \in E(G)$ and $u_j \in N_G(x) \cap V_j$.
Since $N_G(x) \cap N_G(y) \cap V_i$ is anti-complete to $N_G(x) \cap V_j$, we see that $u_i \in (N_G(y) \less N_G(x)) \cap V_i$.
This proves Lemma~\ref{lem:KempeNbrs}(a).\medskip

To prove Lemma~\ref{lem:KempeNbrs}(b), let $H = \overline{G[N_G(x)]}$.
Assume that $d_G(x) = t + 1$ and that $y$ belongs to a cycle, say $C_k$, of $H$, where $k \ge 5$.
By Proposition~\ref{prop:DCBasics}(j), $d_G(y) \ge t + 2$, and by Proposition~\ref{prop:DCBasics}(i), $H$ is the union of isolated vertices and cycles of length at least five.
Clearly, $|N_G(x) \cap N_G(y)| = t - 2$.
By Proposition~\ref{prop:DCBasics}(c), we may assume that $V_i \cap (N_G(x) \cap N_G(y)) = \{v_i\}$ for all $i \in \{1, \dots, t - 2\}$.
Then $N_G(x) \cap N_G(y) = \{v_1, \dots, v_{t - 2}\}$.
Let $\{a, b\} = N_G(x) \less N_G[y]$.
Since $a$ and $b$ are neighbors of $y$ in a cycle of length at least $5$ in $H$, $ab \in E(G)$.
We may further assume that $a \in V_1$ and $b \in V_2$.
Then $v_1 a y b v_2$ forms a path on five vertices of $C_k$, since $v_1, a \in V_1$ and $v_2, b \in V_2$.
If $k \ge 6$, then $v_1 v_2 \in E(G)$ and both $v_1$ and $v_2$ have precisely one non-neighbor in $\{v_3, v_4, \dots, v_{t - 2}\}$.
We may assume that $v_1 v_3 \notin E(G)$ and $v_2 v_{\ell} \notin E(G)$, where $\ell = 3$ if $k = 6$, and $\ell = 4$ if $k \ge 7$.
For any $i, j \in \{3, 4, \dots, t - 2\}$ with $i \ne j$, if $v_i v_j \notin E(G)$, then by Lemma~\ref{lem:KempeNbrs}(a), there exists $v_j' \in V_j \less v_j$ such that $v_j' y \in E(G)$.
By symmetry, there exists $v_i' \in V_i \less v_i$ such that $v_i' y \in E(G)$.
Therefore, if $C$ is any cycle in $H \setminus V(C_k)$ and $V_m \cap V(C) \ne \emptyset$ for some $m \in \{3, 4, \dots, t - 2\}$, then $y$ is adjacent to a vertex in $V_m \setminus v_m$.
\medskip
 
Assume that $k = 5$. Then $v_1v_2\notin E(G)$ and so  $d_G(y) \ge |N_G(x) \cap N_G(y)| + |\{x\}| + e(H \less V(C_k)) = (t - 2) + 1 + (e(H) - 5) = t + e(H) - 6$.
Next assume that $k = 6$.
Then $v_{\ell} = v_3$.
Since both $N_G(x) \cap N_G(y) \cap V_1$ and $N_G(x) \cap N_G(y) \cap V_2$ are anti-complete to $N_G(x) \cap V_3$, by Lemma~\ref{lem:KempeNbrs}(a), $N_G(y)  \cap (V_1\less\{a, v_1\})\ne\emptyset$  and $N_G(y)  \cap (V_2\less\{b, v_2\})\ne\emptyset$.  
 Then $d_G(y)  \ge  |N_G(x) \cap N_G(y)| + |\{x\}| +|N_G(y)  \cap (V_1\less\{a, v_1\})|+ |N_G(y)  \cap (V_2\less\{b, v_2\})|+e(H \less V(C_k)) \ge (t - 2) + 1+1+1 + (e(H) - 6) = t + e(H) - 5$.
Finally assume that $k \ge 7$.
Then $v_{\ell} = v_4$.
Since  $N_G(x) \cap N_G(y) \cap V_1$ is  anti-complete to $N_G(x) \cap V_3$ and $N_G(x) \cap N_G(y) \cap V_2$ is  anti-complete to $N_G(x) \cap V_4$, by Lemma~\ref{lem:KempeNbrs}(a), $N_G(y)  \cap (V_1\less\{a, v_1\})\ne\emptyset$  and $N_G(y)  \cap (V_2\less\{b, v_2\})\ne\emptyset$.     
 As observed earlier,  for any $i, j \in \{3, 4, \dots, t - 2\}$ with $i \ne j$ and $v_i v_j \notin E(G)$, $y$ has at least one neighbor in each of $V_i \less v_i$ and $V_j \less v_j$ in $G$. 
 Hence 
 $d_G(y)  \ge  |N_G(x) \cap N_G(y)| + |\{x\}|+|N_G(y)  \cap (V_1\less\{a, v_1\})|+ |N_G(y)  \cap (V_2\less\{b, v_2\})| +|V(C_k)\less\{a, b, v_1, v_2, y\}|+ e(H \less V(C_k)) \ge (t - 2) + 1+1+1 +(k-5)+ (e(H) - k) = t + e(H) - 4$.
Note that since $k \ge 7$, we see that $e(H) \ge 7$, and so $d(y) \ge t + e(H) - 4> t + 2$.
This completes the proof of Lemma~\ref{lem:KempeNbrs}(b).
\proofsquare
\medskip

\begin{lem}\label{lem:Delta}
Let $G$ be a double-critical, $t$-chromatic graph with $t \ge 6$.
If $G$ is claw-free, then for any $x \in V(G)$, $d_G(x) \le 2t - 4$.
Furthermore, if $d_G(x) < |V(G)| - 1$, then $d_G(x) \le 2t - 6$.
\end{lem}

\pf  Let $x \in V(G)$ be a vertex of maximum degree in $G$, and let $uv$ be any edge of $G \less x$.
Let $c$ be any $(t - 2)$-coloring of $G \less \{u, v\}$ with color classes $V_1, V_2, \dots, V_{t - 2}$.
We may assume that $x \in V_{t - 2}$.
Since $G$ is claw-free, $x$ can have at most two neighbors in each of $V_1, \dots, V_{t - 3}$.
Additionally, $x$ may be adjacent to $u$ and $v$ in $G$.
Therefore $d_G(x) \le 2t - 4$.
If $d_G(x) < |V(G)| - 1$, then $\chi(G[N_G(x)]) \le t - 3$ by Proposition~\ref{prop:DCBasics}(h).
Since $G$ is claw-free, each color class in any $(t - 3)$-coloring of $G[N_G(x)]$ can contain at most two vertices, and so $d_G(x) \le 2t - 6$.
\proofsquare
\bigskip

It is now an easy consequence of Proposition~\ref{prop:DCBasics} and Lemma~\ref{lem:Delta} that Theorem~\ref{thm:DCClawFree678} is true for $t = 6, 7$.
\medskip

\noindent{\bf Proof of Theorem~\ref{thm:DCClawFree678} for $t = 6, 7$.}
Let $G$ and $t \in \{6, 7\}$ be as given in the statement.
Suppose that $G \ne K_t$.
By Proposition~\ref{prop:DCBasics}(d), there exists an edge $xy \in E(G)$ such that $x$ and $y$ have a common non-neighbor.
By Proposition~\ref{prop:DCBasics}(f) and Lemma~\ref{lem:Delta}, $t + 1 \le d_G(x) \le 2t - 6$ and $t + 1 \le d_G(y) \le 2t - 6$.
Thus $t=7$ and $d_G(x) = d_G(y) = 8$, which contradicts Proposition~\ref{prop:DCBasics}(j).
\proofsquare

\section{Proofs of Main Results}\label{sec:DCClawFree678}

In this section, we prove our main results, namely, Theorem~\ref{thm:t+1nbrs} and Theorem~\ref{thm:DCClawFree678} for the case $t = 8$.
We first prove Theorem~\ref{thm:t+1nbrs}.
\medskip

\noindent{\bf Proof of Theorem~\ref{thm:t+1nbrs}.}
Let $G$ and $x$ be as given in the statement.
Let $H = \overline{G[N_G(x)]}$.
Then $|V(H)| = t + 1$.
Note that if $d_G(x)=|V(G)|-1$, then it follows from Proposition~\ref{prop:DCBasics}(f) that $G = K_{t + 1}$, a contradiction. Thus $d_G(x)<|V(G)|-1$. 
Now by Proposition~\ref{prop:DCBasics}(g) and Proposition~\ref{prop:DCBasics}(h) applied to the vertex $x$, $\alpha(\overline{H}) \le 2$ and $\chi(\overline{H}) \le t - 3$.
Let $c^*$ be any $(t - 3)$-coloring of $\overline{H}$.
Then each color class of $c^*$ contains at most two vertices.
Since $|V(H)| = t + 1$, we see that at least four color classes of $c^*$ must each contain two vertices.
By Proposition~\ref{prop:DCBasics}(e),  $H$ has  at least eight vertices of degree two  and so $e(H) \ge 8$.
This proves Theorem~\ref{thm:t+1nbrs}(a).
\medskip

To prove Theorem~\ref{thm:t+1nbrs}(b), let $y \in N_G(x)$.
Since $d_G(x) = t + 1$, by Proposition~\ref{prop:DCBasics}(i), either $|N_G(x) \cap N_G(y)| = t$ or $|N_G(x) \cap N_G(y)| = t - 2$.
Assume that $|N_G(x) \cap N_G(y)| = t - 2$. Then $y$ belongs to a cycle of length $k\ge5$  in $H$ because $H$ is a disjoint union of isolated vertices and cycles. 
By Proposition~\ref{prop:DCBasics}(i), $y$ belongs to a cycle of length at least $5$ in $H$.
By Theorem~\ref{thm:t+1nbrs}(a), $e(H) \ge 8$.
Note that if  $5 \le k \le 7$, then by Proposition~\ref{prop:DCBasics}(i), $H$ has at least two cycles of length at least 5, and so $e(H) \ge k + 5 \ge 10$.
Thus by Lemma~\ref{lem:KempeNbrs}(b), $d_G(y) \ge t + 4$.
If $d_G(y) = t + 4$, then it follows from Lemma~\ref{lem:KempeNbrs}(b) that  either $k=8$  and $H$ is isomorphic to $C_{8} \cup \overline{K}_{t-7}$  or   $k=5$  and $H$ is isomorphic to $C_{5} \cup C_5\cup \overline{K}_{t-9}$.  
So we may assume that $|N_G(x) \cap N_G(y)| = t$.
Let $c$ be any $(t - 2)$-coloring of $G \less \{x, y\}$ with color classes $V_1, V_2, \dots, V_{t - 2}$.
Since $\alpha(\overline{H}) \le 2$, we may further assume that $N_G(x) \cap V_1=\{v_1, v_1'\}$, $N_G(x) \cap V_2 = \{v_2,v_2'\}$,  and $N_G(x) \cap V_i = \{v_i\}$ for all $i \in \{3, 4, \dots, t - 2\}$.
Then $v_1v_1', v_2v_2'\in E(H)$. By Proposition~\ref{prop:DCBasics}(i) applied to the vertex $x$,  $e_H(\{v_1, v_1', v_2, v_2'\}, \{v_3, v_4, \dots, v_{t-2}\})\le4$.  
By Theorem~\ref{thm:t+1nbrs}(a), $e(H) \ge 8$.  Thus there must exist at least four vertices in $\{v_3, v_4, \dots, v_{t - 2}\}$, say $v_3, v_4, v_5, v_6$, such that $d_H(v_j)=2$ and $y$ is  adjacent to at least one vertex of $V_j \less v_j$ in $G$ for all $j\in\{3,4,5,6\}$.  
Therefore $|N_G(y) \less N_G[x]| \ge 4$ and so $d_G(y) = |N_G[x] \cap N_G(y)| + |N_G(y) \less N_G[x]|\ge (t + 1) + 4 = t + 5$.

\medskip

This completes the proof of Theorem~\ref{thm:t+1nbrs}.
\proofsquare
\bigskip


We are now ready to complete the proof of Theorem~\ref{thm:DCClawFree678}.
\medskip

\noindent{\bf Proof of Theorem~\ref{thm:DCClawFree678} for $t = 8$.}
Let $G$ and $t = 8$ be as given in the statement.
Suppose that $G \ne K_8$.
We claim that
\medskip

\noindent {\bf Claim \refstepcounter{counter}\label{e:10Reg}  \arabic{counter}.} $G$ is 10-regular.

\pf  By Lemma~\ref{lem:DomVertex} and Theorem~\ref{thm:DCClawFree678} for $t = 7$, $\Delta(G) \le |V(G)| - 2$.
By Proposition~\ref{prop:DCBasics}(f) and Lemma~\ref{lem:Delta}, we see that $9 \le d_G(x) \le 10$ for all vertices $x \in V(G)$.
By Corollary~\ref{cor:t+1nbrs}, $G$ is $10$-regular.
\proofsquare
\medskip

\noindent {\bf Claim\refstepcounter{counter}\label{e:2NonNbrs}  \arabic{counter}.} For any $x \in V(G)$, $2 \le \delta(\overline{G[N_G(x)]}) \le \Delta(\overline{G[N_G(x)]}) \le 3$.

\pf Let $x \in V(G)$.
Then $x$ has at least one non-neighbor in $G$, otherwise $G = K_{11}$ by Claim~\ref{e:10Reg}, a contradiction.
By Proposition~\ref{prop:DCBasics}(h), $\chi(G[N_G(x)]) \le 5$.
Since $G$ is claw-free, we see that $\alpha(G[N_G(x)]) = 2$, and so $\chi(G[N_G(x)]) = 5$ since every color class can contain at most two vertices.
Thus every vertex of $N_G(x)$ has at least one non-neighbor in $G[N_G(x)]$.
By Proposition~\ref{prop:DCBasics}(e) and Proposition~\ref{prop:DCBasics}(c), $2 \le \delta(\overline{G[N_G(x)]}) \le \Delta(\overline{G[N_G(x)]}) \le 3$.
\proofsquare
\bigskip

\noindent {\bf Claim\refstepcounter{counter}\label{e:Not2Reg}  \arabic{counter}.} For any $x\in V(G)$, $\Delta(\overline{G[N_G(x)]}) = 3$.
That is, $\overline{G[N_G(x)]}$ is not 2-regular.

\pf Suppose that there exists a vertex $x \in V(G)$ such that $\overline{G[N_G(x)]}$ is $2$-regular.
Let $y \in N_G(x)$ and let $c$ be any $6$-coloring of $G \less \{x, y\}$ with color classes $V_1, V_2, \dots, V_6$.
Let $W = N_G(x) \cap N_G(y)$.
Then $|W| = 7$ because $\overline{G[N_G(x)]}$ is $2$-regular.
By Proposition 2.1(c), we may assume that $|V_1 \cap W| = 2$ and $|V_i \cap W| = 1$ for every $i \in \{2, 3, 4, 5, 6\}$.
Let $V_1 \cap W = \{v_1, u_1\}$ and $V_i \cap W = \{v_i\}$ for each $i \in \{2, 3, 4, 5, 6\}$.
Since $G$ is claw-free, we may further assume that $N_G(x) \cap V_2 = \{v_2, u_2\}$ and $N_G(x) \cap V_3 = \{v_3, u_3\}$.
Clearly, $y u_2, y u_3 \notin E(G)$ and thus $u_2 u_3 \in E(G)$ because $G$ is claw-free.
Since $\overline{G[N_G(x)]}$ is $2$-regular, we see that $G[\{v_4, v_5, v_6\}]$ is not a clique.
We may assume that $v_4 v_5 \notin E(G)$.
By Lemma~\ref{lem:KempeNbrs}(a), $N_G(y) \cap (V_j \less \{v_j\}) \ne \emptyset$ for all $j \in \{4, 5\}$.
Let $w_4 \in V_4 \less v_4$ and $w_5 \in V_5 \less v_5$ be two other neighbors of $y$ in $G$.
Then $N_G(y) \less N_G[x] = \{w_4, w_5\}$ since $G$ is 10-regular by Claim~\ref{e:10Reg}.
By Lemma~\ref{lem:KempeNbrs}(a), $v_6$ must be complete to $\{v_2, v_3, v_4, v_5\}$ in $G$.
Notice that $v_6$ is complete to $\{u_2, u_3\}$ in $G$ since $\overline{G[N_G(x)]}$ is 2-regular.
Thus $v_6$ must be anti-complete to $\{v_1, u_1\}$ in $G$ and so $G[\{x, v_1, u_1, v_6\}]$ is a claw, a contradiction.
\proofsquare
\bigskip

From now on, we fix an arbitrary vertex $x\in V(G)$.
Let $H = \overline{G[N_G(x)]}$.
By Claim~\ref{e:Not2Reg}, let $y \in N_G(x)$ with $|N_G(x) \cap N_G(y)| = 6$.
We choose such a vertex $y \in N_G(x)$ so that $N_G(x) \less N_G[y]$ contains as many vertices of degree two in $H$ as possible.
Let $c$ be any $6$-coloring of $G \less\{x, y\}$ with color classes $V_1, V_2,\dots, V_6$.
We may assume that $V_i \cap N_G(x) \cap N_G(y) = \{v_i\}$ for all $i \in \{1, 2, 3, 4, 5, 6\}$.
Since $G$ is claw-free, we may further assume that $N_G(x) \cap V_j = \{v_j, u_j\}$ for all $j \in \{1, 2, 3\}$.
Notice that $y$ is anti-complete to $\{u_1, u_2, u_3\}$ in $G$,  and since $G$ is claw-free, $G[\{u_1, u_2, u_3\}] = K_3$.
Let $A = \{u_1, u_2, u_3\}$, $B = \{v_1, v_2, v_3\}$, and $C = \{v_4, v_5, v_6\}$.
\bigskip

\noindent {\bf Claim\refstepcounter{counter}\label{e:BNotCompC}  \arabic{counter}.} $B$ is not complete to $C$ in $G$.

\pf Suppose that $B$ is complete to $C$ in $G$.
Then $e_H(C, A) = \sum_{v \in C} d_H(v) - 2e(H[C]) \ge 6 - 2 e(H[C])$.
For each $i \in \{1, 2, 3\}$, $u_i v_i, u_i y \notin E(G)$ and $d_H(u_i) \le 3$.
Thus $e_{H}(A, C) \le 3$ and so $e(H[C]) \ge 2$.
Since $G$ is claw-free, we have $e(H[C]) =2$.
We may assume that $v_4 v_6 \notin E(H)$.
Then $v_4 v_6 \in E(G)$ and $v_4 v_5, v_5 v_6 \notin E(G)$.
Since $d_H(v_4) \ge 2$,   $d_H(v_6) \ge 2$,  and $B$ is complete to $C$ in $G$, we may assume that $u_2 v_4, u_3 v_6 \notin E(G)$.
Note that $H$ is not 3-regular since $e_H(A, C) \le 3$ and $e_H(B, C) = 0$.
By the choice of $y$, $d_H(u_1) = 2$ and $d_H(v_j) = 2$ for all $j \in \{4, 5, 6\}$.
Since $d_H(u_2) = d_H(u_3) = 3$, by the choice of $y$ again, $d_H(v_2) = d_H(v_3) = 3$.
Thus $G[B] = \overline{K_3}$ and so $G[\{x\} \cup B]$ is a claw, a contradiction.
\proofsquare
\bigskip

\noindent {\bf Claim \refstepcounter{counter}\label{e:CComp} \arabic{counter}.} $G[C] = K_3$.

\pf Suppose that $G[C]$ contains a missing edge, say $v_4 v_5 \notin E(G)$.
By Lemma~\ref{lem:KempeNbrs}(a), there exist $w_4 \in V_4 \less v_4$ and $w_5 \in V_5 \less v_5$ such that $y w_4, y w_5 \in E(G)$.
By Claim~\ref{e:BNotCompC}, we may assume that $v_3 v_j \notin E(G)$ for some $j \in \{4, 5, 6\}$.
By Lemma~\ref{lem:KempeNbrs}(a), $y$ has another neighbor, say $w_3$, in $V_3 \less v_3$.
Since $G$ is 10-regular, $\{w_3, w_4, w_5\} = N(y) \less N[x]$, so by Lemma~\ref{lem:KempeNbrs}(a), $v_4 v_5$ is the only missing edge in $G[C]$ and $\{v_1, v_2\}$ is complete to $C$ in $G$.
If $e_H(A, C) = 3$, then $d_H(u_i) = 3$ for all $i \in \{1, 2, 3\}$.
By the choice of $y$, $d_H(v_3) = 3$, or else we could replace $y$ with $u_3$.
Notice that for all $i \in \{4, 5, 6\}$, $e_H(v_i, A \cup \{v_3\}) \ge 1$, and so by the choice of $y$, $d_H(v_i) = 3$, or else we could replace $y$ with $v_3$.
Thus $e_H(A, C) \ge 5$, which is impossible.
Hence  $e_H(A, C) \le 2$. Notice that $e_H(A, C) = (d_H(v_4) - 1) + (d_H(v_5) - 1) + d_H(v_6) - e_H(v_3, C) \ge 2$.  it follows  that $e_H(A, C) = 2$, $e_H(v_3, C) = 2$ and $d_H(v_i) = 2$ for all $i \in \{4, 5, 6\}$.
Then $N_G(x) \less N_G[y]$ has at most one vertex of degree two in $H$, but $N_G(x) \less N_G[v_3]$ has two vertices of degree two in $H$, contradicting the choice of $y$.
\proofsquare
\bigskip

\noindent {\bf Claim\refstepcounter{counter}\label{e:ABedges}  \arabic{counter}.} $v_1 u_1, v_2 u_2$, and $v_3 u_3$ are the only edges in $H[A \cup B]$.

\pf Suppose that $H[A \cup B]$ has at least four edges.
By Claim~\ref{e:CComp} and Claim~\ref{e:2NonNbrs}, $e_H(A \cup B, C) \ge 6$.
On the other hand, $e_H(A \cup B, C) = \sum_{v \in A \cup B} d_H(v) - 2e(H[A \cup B]) - 3 \le 15 - 2e(H[A \cup B])$.
It follows that $e(H[A \cup B]) = 4$ and $A \cup B$ contains at most one vertex of degree two in $H$.
Thus $e_H(A \cup B, C) \le 7$ and so at least two vertices of $C$, say $v_4$ and $v_5$,  are of degree two in $H$.
Since $e_H(A, C) \le 3$ and $G[C] = K_3$ by Claim~\ref{e:CComp}, we may assume that $v_4 v_3 \notin E(G)$.
If $d_H(v_3) = 3$, then since $d_H(v_4) = 2$ and at most one vertex of $A \cup B$ has degree two in $H$, by the choice of $y$, exactly one of $u_1, u_2, u_3$ has degree two in $H$.
Then $e_H(A \cup B, C) = 6$.
Thus $d_H(v_j) = 2$ for all $j \in \{4, 5, 6\}$ and by the choice of $y$, each vertex of  $B$ is adjacent to at most one vertex of $C$ in $H$.
Thus $e_H(A \cup B, C) \le 5$, a contradiction.
Hence $d_H(v_3) = 2$.
Now $d_H(u_i) = 3$ for all $i \in \{1, 2, 3\}$ because at most one vertex of $A \cup B$ has degree two in $H$.
We see that $N(x) \less N[y]$ has no vertex of degree two in $H$ but $N(x) \less N[u_3]$ has at least one vertex of degree two in $H$, contrary to the choice of $y$.
\proofsquare
\bigskip

By Claim~\ref{e:ABedges}, we see that for any $i \in \{1, 2, 3\}$, $v_i v_j \notin E(G)$ for some $j \in \{4, 5, 6\}$.
By Lemma~\ref{lem:KempeNbrs}(a), let $w_i \in V_i \less v_i$ be such that $y w_i \in E(G)$ for all $i \in \{1, 2, 3\}$.
Let $D = \{w_1, w_2, w_3\}$.
Then $N_G(y) \less N_G[x] = D$ and $G[D] = K_3$ because $G$ is claw-free.
Clearly, $D$ is not complete to $C$ in $G$, otherwise $G[\{y\} \cup D \cup C] = K_7$, contrary to Proposition~\ref{prop:DCBasics}(a).
We may assume that $w_3 v_4 \notin E(G)$.
For each $i \in \{1, 2\}$, $v_i v_3, v_i u_3 \in E(G)$ by Claim~\ref{e:ABedges}.
Thus $v_1 w_3, v_2 w_3 \notin E(G)$ because $G$ is claw-free.
Notice that $w_3, x, v_1, v_2, v_4 \in N_G(y)$ and $w_3$ is anti-complete to $\{x, v_1, v_2, v_4\}$ in $G$.
Thus $\Delta(\overline{G[N_G(y)]}) \ge 4$, contrary to Claim~\ref{e:2NonNbrs}.

This completes the proof of Theorem~\ref{thm:DCClawFree678}.
\proofsquare
\bigskip

\section*{Acknowledgement}
The authors would like to thank the anonymous referees for  many helpful comments, which  greatly improve the presentation of this paper.

\end{document}